\documentclass{article}
\usepackage{float}
\usepackage{amsmath,amssymb,amsfonts,textcomp,amsthm}
\usepackage{epsfig,epstopdf}
\usepackage{mathrsfs}
\usepackage{caption}
\usepackage{graphicx}

\newcommand{\parcial}[2]{\frac{\partial #1}{\partial #2}}
\newcommand{\norm}[1]{\left\lVert#1\right\rVert}

\title{Random noise attenuation on finite-difference wave propagation using fuzzy transform}

\author{Rafael Aleixo, UFSC }

\begin{document}
\newtheorem{defin}{Definition}[section]
\newtheorem{teo}{Theorem }
\newtheorem{coro}{Corollary}
\newtheorem{prop}{Proposition }
\newtheorem{lem}{Lema }
\newtheorem{ejem}{Example}
\newtheorem{nota}{Note}
\def\demo{\textit{Proof}\quad}
\def\udots{\rotatebox{90}{$\ddots$}}
\pdfinclusioncopyfonts=1

\maketitle

\begin{abstract}

Fuzzy Transform (F-transform) has been introduced as an approximation method which encompasses both classical transforms and approximation methods studied in fuzzy modeling and fuzzy control.  It has been proved that, under some conditions, F-transform can remove a periodical noise and it can significantly reduce random noise. In this work we apply the F-transform methodology to propose a finite-difference approach to solving the acoustic wave equation with noisy initial conditions.

\end{abstract}

\section{Introduction}

Fuzzy Transform (F-transform) has been introduced as an approximation method which encompasses both classical transforms as well as approximation methods studied in fuzzy modelling and fuzzy control \cite{Perfilieva1,Perfilieva3,Perfilieva2}. Since then, many valuable properties of this method has been shown and as consequence many applications has been publishing from geology \cite{Perfilieva5} to electrical engennering \cite{Stepnicka3}, including image processing \cite{Perfilieva2,Stepnicka2} and differential equations \cite{Stepnicka5,Stepnicka,Holcapek:2017}. As in other approximations technics, it can be shown that if the original function is replaced by an approximation model, then a certain simplification of complex computations could be achieved but additionally F-transform has been demonstrated to be a robust method. It has been proved that, under some conditions, it can significantly reduce random/periodic noise \cite{Perfilieva4,Stepnicka}.

Following the work started by Aleixo and Amazonas \cite{Aleixo:2016} we investigate in more details the nature of the F-transform when applied on wave propagations. The main idea consists in the replacement of an continuous function on a real closed interval by its discrete representation (using the direct F-transform). According to specific situation, computations or calculus are made on discrete representation of function and then transformed back to the space of continuous functions (using the inverse F-transform). The result obtained by applying both F-transforms is a good simplified approximation of original function. The main application of the ideas developed in this work is on reverse time migration, where the back-propagated data presents, in general, random and periodic noise contents. Besides, these ideas can be applies to multiple attenuation algorithms.

\section{The Fuzzy Transform}

In this section we introduce the definition of F-transform and its inverse. Consider a real closed interval $[a,b]$ as a universe. First of all, we introduce the definition of basic functions, that play a key role in the definition of F-transform.

Let $x_{i}=a+h(i-1)$ be nodes (a uniform partition) on $[a, b]$ where $h = (b - a)/(n - 1)$, $n \geq 2$ and $i = 1,\ldots,n$. Let $A_1(x), . . . , A_n(x)$ be real-valued functions defined on $[a, b]$. These functions are called basic functions if each of them satisfies the following conditions:
\begin{itemize}
\item $\displaystyle A_i : [a, b] \to [0, 1]$ and $A_i(x_i) = 1;$
\item $\displaystyle A_i(x) = 0$ if $x \notin (x_{i-1}, x_{i+1})$, where $x_0 = a, x_{n+1} = b;$
\item $\displaystyle A_i(x)$ is continuous, strictly increases on $[x_{i-1}, x_i]$ and strictly decreases on $[x_i, x_{i+1}]$;
\item $\displaystyle\sum_{i=1}^{n}A_i(x_i)=1$, for all $x \in [a, b]$;
\item $\displaystyle A_i(x_i-x) = A_i(x_i+x)$, for all $x \in [0, h]$, $i = 2, \ldots , n- 1,$
\item $\displaystyle A_{i+1}(x) = A_i(x - h)$, for all $x \in [a + h, b], i =
2, \ldots , n - 2$.
\end{itemize}

The partition $x_{i}=a+h(i-1)$ and the basic functions $A_1(x), . . . , A_n(x)$ are called an uniform fuzzy partition. A typical example of basic functions are uniform triangular functions. In some cases, extended fuzzy partitions might be considered \cite{Stefanini:2011}. Now, the notion of F-transform and its inverse in two variables will be defined, the generalization to three or more variables and the one-dimensional version is straightforward.

Let $f(x, y)$ be an arbitrary real-valued continuous function on $\mathcal{D} = [a, b] \times [c, d]$. Let $A_1(x),\ldots,A_n(x)$ on $[a, b]$ and $B_1(y),\ldots,B_m(y)$ on $[c, d]$ form uniform fuzzy partitions, not necessarily the same. We say that a matrix $[F_{ij} ]_{nm}$ of real numbers is the F-transform of $f(x, y)$ with respect to the given basic functions if
$$
F_{ij}=\frac{\displaystyle\int_{c}^{d}\int_{a}^{b} f(x,y)A_i(x)B_j(y)dxdy}{\displaystyle\int_{c}^{d}\int_{a}^{b} A_i(x)B_j(y)dxdy}, \,\,\,i= 1, \ldots, n,\,j= 1, \ldots , m.
$$
Moreover, let $[F_{ij}]_{nm}$ be the F-transform of a function $f$. Then the function
$$
f^{F}_{n,m}(x,y)=\sum_{i=1}^{n}\sum_{j=1}^{m} A_i(x)B_j(y)F_{ij}
$$
is called the inverse F-transform \cite{Stepnicka}.

An approximation property of the inverse F-transform with a given arbitrary precision is shown by Perfilieva and Val{\' a}{\v s}ek \cite{Perfilieva4}. The statement is: let $f$ a continuous function on $[a,b]$. Then, for any $\epsilon>0$, there exists $n_\epsilon$ and a fuzzy partition $A_1, \ldots , A_{n_\epsilon}$ of $[a,b]$ such that for all $x\in[a,b]$
$$
|f(x)-f^{F}_{n_\epsilon}(x)|<\epsilon,
$$
where $f^{F}_{n_\epsilon}(x)$ is the inverse F-transform of $f$ with respect to the fuzzy partitions $A_1,\ldots , A_{n_\epsilon}$. This property naively says the inverse F-transform is an approximation of $f$ with a given desired precision.

\section{Fuzzy transform and noise attenuation}

The purpose of this work is to show how F-transform impacts finite-difference (FD) solutions of a certain class of partial differential equations (PDEs) when the initial condition is contaminated by noise. Fuzzy transform works well on noise removal. An extensive work was done to identify which kind of noise can be removed by using F-transforms \cite{Perfilieva4}.

\begin{figure}[H]
\centerline{ \includegraphics[trim={0cm 0cm 0cm 1cm},clip,scale=0.6]{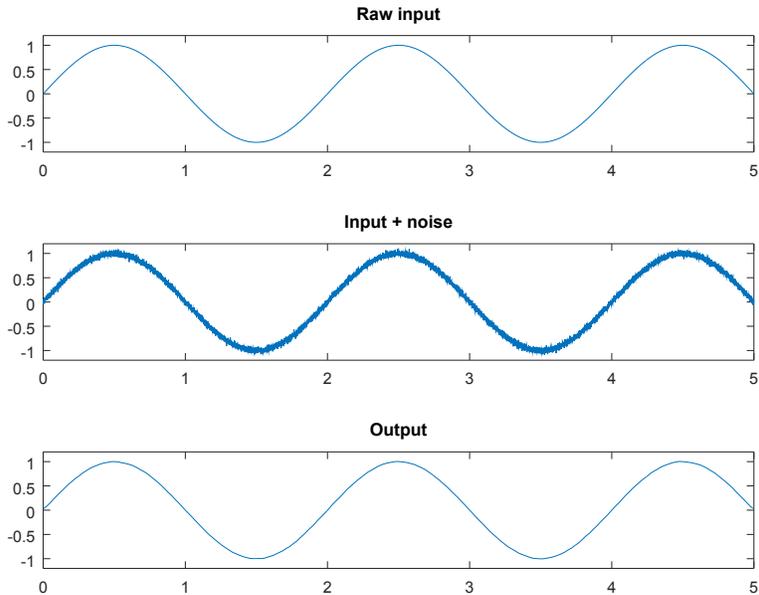}}
\caption{F-transform methodology validation.}
\label{fig0}
\end{figure}

Perfilieva and Val{\' a}{\v s}ek \cite{Perfilieva4} showed the application of direct and then inverse Fuzzy transform to a noisy function works as a denoise procedure. Among the variety of different noises, this methodology can be applied to random noise attenuation. Figure \ref{fig0} shows the effect of the Fuzzy transform as a tool for random noise attenuation. Firstly, a raw sine function (top graph) is presented, as reference. Then, random noise is added to the function (mid graph). After that, the F-transform methodology is applied to the noise-contaminated sine function (mid graph). It consists on an application of a forward F-transform along with an inverse F-transform. The bottom graph shows the result.

The methodology presented above is not what exactly what we intend to do, but it validates the idea we want to explore. In fact, the methodology we propose for solving PDEs numerically is a variation of the forward+inverse F-transform idea. We construct the forward F-transform of the PDE's solution and then we apply the inverse F-transform. Next section we present how the forward F-transform is accomplished.

\section{Fuzzy transform and finite-difference schemes}

F-transform can be used to solve PDEs numerically \cite{Stepnicka,Holcapek:2017}, specifically by finite-difference (FD) schemes. Moreover, F-transform approach, on FD schemes, works well when noise is presented on the external force factor \cite{Stepnicka}. In this work our approach is slightly different, we assume the initial condition is perturbed by noise.

For the rest of the text we assume the considered functions have as many continuous derivatives as necessary to the calculations below to be valid. In general, the Fuzzy transform can be used on the calculation of FD solutions for the class of differential equations
\begin{equation}\label{eq}
\sum_{n=1}^{N}\alpha_n\frac{\partial{^n}}{\partial t^n}u+\sum_{m=1}^{M}\beta_m\frac{\partial{^m}}{\partial x^m}f(u)=q(x,t),
\end{equation}
where $N$, $M$, $\alpha_n$ and $\beta_n$ are finite constants, $q$ is the external force and $f$ is a continuous function. Special cases of equation (\ref{eq}) are the hyperbolic conservation laws, acoustic wave equation, Burgers' equation, etc. As the acoustic wave equation is widely studied in geophysics, we present how to use F-transform on this equation. The calculations for the general case are quite tedious. The ideas presented for the wave equation are enough to understand the process. Besides, for the sake of notation, the process is applied to a fixed FD scheme. The general case is straightforward because of the F-transform linearity. We present the theory for one-dimensional wave equation because the general description is quite fastidious.

The derivation described below follows the ideas presented in \cite{Stepnicka}. Let us consider a set $D=[a,b]\times[0,T]$  and let $u(x,t)$ be a solution of the wave equation on $D$,
\begin{equation}
\parcial{^2 u}{x^2} - \frac{1}{c^2}\parcial{^2 u}{t^2}=q(x,t),
\label{wave}
\end{equation}
with the following initial and boundary conditions
$$
\begin{array}{l}
  u(x,0)=f(x), \qquad \displaystyle\parcial{u}{t}(x,0)=g(x),\\\\
  \qquad u(a,t)=T_1(t), \qquad u(b,t)=T_2(t).
\end{array}
$$

After applying the F-transform, equation (\ref{wave}) is turned into the following algebraic equation
\begin{equation}
\label{recur1}
U^{xx}_{ij}-\frac{1}{c^2} U^{tt}_{ij}=Q_{ij}
\end{equation}
where $U^{xx}_{ij}, U^{tt}_{ij}$ and $Q_{ij}$ are the F-transform components of the functions $\displaystyle \parcial{^2 u}{x^2}, \parcial{^2 u}{t^2}$ and $q$ respectively.

The ideia is to use Taylor-series expressions to estimate the F-transforms of $U^{xx}_{ij}, U^{tt}_{ij}$. Using finite differences to approximate the partial derivatives on left side of (\ref{wave}) we would have
\[\parcial{^2 u}{x^2}=\frac{u(x+h_x,t)-2u(x,t)+u(x-h_x,t)}{h^2_x} + O(h^2_x)\]
and
\[\parcial{^2 u}{t^2}=\frac{u(x,t+h_t)-2u(x,t)+u(x,t-h_t)}{h^2_t} + O(h^2_t).\]

Now, for $j=2,\ldots,m-2$, we can consider, without loss of generality, the function $u(x,t+h_t)A_i(x)B_j(t)$ being defined on the entire interval $[0,T]$, so we have that

$$
\begin{array}{l}
  \displaystyle\int_a^b \int_0^T u(x,t+h_t)A_i(x)B_j(t)dtdx=\\\\
  \displaystyle\qquad =\int_a^b\int_{t_j}^{t_{j+2}} u(x,t)A_i(x)B_j(t-h_t)dtdx \\\\
  \displaystyle\qquad =\int_a^b\int_{t_j}^{t_{j+2}} u(x,t)A_i(x)B_{j+1}(t)dtdx.
\end{array}
$$
Therefore,
$$
\begin{array}{l}
  \displaystyle\int_a^b \int_0^T u(x,t+h_t)A_i(x)B_j(t)dtdx=\\\\
  \displaystyle\qquad =\int_a^b\int_{0}^{T} u(x,t)A_i(x)B_{j+1}(t)dtdx.
\end{array}
$$

This last expression it is valid also when $j=1$ or $j=m-1$, because definition of uniform fuzzy partition. It is also easy to prove that
$$
\begin{array}{l}
\displaystyle\int_a^b \int_0^T u(x,t-h_t)A_i(x)B_j(t)dtdx\\\\
\displaystyle\qquad=\int_a^b\int_{0}^{T} u(x,t)A_i(x)B_{j-1}(t)dtdx,
\end{array}
$$
for $j=1,\ldots,m-1$. In the same way we have
$$\begin{aligned}
\displaystyle \int_a^b \int_0^T A_i(x)B_j(t)dtdx&=\displaystyle \int_a^b\int_{x_{j-1}}^{x_{j+1}}A_i(x)B_j(t)dtdx\\\\
&\displaystyle =\int_a^b\int_{x_j}^{x_{j+2}} A_i(x)B_{j+1}(t)dtdx  \\\\
&\displaystyle =\int_a^b\int_{0}^{T} A_i(x)B_{j+1}(t)dtdx.
\end{aligned}$$
and
$$
\displaystyle \int_a^b \int_0^T A_i(x)B_j(t)dtdx=\int_a^b\int_{0}^{T} A_i(x)B_{j-1}(t)dtdx.
$$

Now we can use these previous calculations to estimate $U^{tt}_{ij}$ as follows:
$$
\begin{aligned}
U^{tt}_{ij}&=\frac{\displaystyle \int_a^b \int_0^T \parcial{^2 u}{t^2}A_i(x)B_j(t)dtdx}{\displaystyle \int_a^b \int_0^T A_i(x)B_j(t)dtdx}\\\\
&\approx \frac{\displaystyle \frac{1}{h^2_t} \int_0^T\int_a^b (u(x,t+h_t)-2u(x,t)+u(x,t-h_t))A_i(x)B_j(t)dxdt}{\displaystyle \int_a^b \int_0^T A_i(x)B_j(t)dxdt}
\end{aligned}
$$

So, we obtain
\begin{equation}
\label{recur2}
U^{tt}_{ij}\approx\frac{1}{h^2_t}(U_{i(j+1)}- 2U_{ij}+U_{i(j-1)}),
\end{equation}
and similarly
\begin{equation}
\label{recur3}
U^{xx}_{ij}\approx\frac{1}{h^2_x}(U_{(i+1)j}- 2U_{ij}+U_{(i-1)j}).
\end{equation}

Now, joining this approximations and using equations (\ref{recur1}), (\ref{recur2}) and (\ref{recur3}) we can obtain the following recursive equation
\begin{equation}
\label{recur4}
\begin{array}{lll}
  U_{i(j+1)} & = & r^2[U_{(i-1)j}+U_{(i+1)j}]+ \\\\
   & + & 2(1-r^2)U_{ij}-U_{i(j-1)}-c^2h^2_tQ_{ij},
\end{array}
\end{equation}
where $r=\displaystyle\frac{c h_t}{h_x}$ and $i=1,\ldots,n$, $j=1,\ldots,m$.

We still need to deal with initial and boundary conditions, so for example, using initial condition $u(0,t)=T_1(t)$ we have
$$
\begin{aligned}
U_{1j}&=\frac{\displaystyle \int_a^b\int_{0}^{T} u(x,t)A_1(x)B_{j}(t)dtdx}{\displaystyle \int_a^b\int_{0}^{T} A_1(x)B_{j}(t)dtdx}\\\\
&=\frac{\displaystyle \int_{t_{j-1}}^{t_{j+1}}\int_{a}^{x_2} u(x,t)A_1(x)B_{j}(t)dxdt}{\displaystyle \int_{t_{j-1}}^{t_{j+1}}\int_a^{x_2} A_1(x)B_{j}(t)dxdt}\\\\
&\approx \frac{\displaystyle \int_{t_{j-1}}^{t_{j+1}}\int_a^{x_2} T_1(t)A_1(x)B_{j}(t)dxdt}{\displaystyle \int_{t_{j-1}}^{t_{j+1}}\int_a^{x_2} A_1(x)B_{j}(t)dxdt}=T_1^j,\end{aligned}\\\\
$$
that is, the boundary vector $U_{1j}$ is the F-transform of $T_1(t)$. Analogously,  we obtain $U_{nj}=T_2^j$ and $U_{i1}=F_i$, where $T_2^j$ is the F-transform of $T_2(t)$ and $F_i$ the F-transform of $f(x)$. Lastly, we can use
\[\parcial{u}{t}\approx \frac{u(x,t+h_t)-u(x,t-h_t)}{2h_t}\]
and initial condition to obtain $$G_i=\frac{U_{i2}-U_{i0}}{2h_t},$$ which let us define $U_{i0}$ for $i=2,\ldots,n-1$ as
$U_{i0}=-2h_tG_i+U_{i2}.$

Now, we have a complete description of the recursive equation for the F-transform $U_{ij}$, i.e., we developed a methodology to construct the F-transform of the wave-equation solution. Bear in mind this recursive equation is similar to the finite-difference recursive equation for analogous partial derivatives, which means the approximation process for the F-transform is equivalent to a FD method, even if they are conceptually different. We can take advantage of this similarity and, in fact, we can say, for example, the numerical stability of algorithm (\ref{recur4}) is achieved if $0<r\leqslant 1$. Moreover, it is easy to prove (\ref{recur4}) converges using numerical analysis technics \cite{Thomas}. After running (\ref{recur4}), the result is the solution of the wave equation at the F-transformed domain. To return to the original domain, an inverse F-transform must be applied. When the $q(x,y)=0$, the F-transform recursive equation and FD recursive equation are exactly the same, except for the initial and boundary conditions. The derivation of the F-transform recursive equations for the class of equations (\ref{eq}) is straightforward, and when the external force is zero, FD and F-transform algorithms are the same. Which means, the F-transform methodology is basically a FD scheme plus an inverse F-transform application. Then, the noise reduction is achieved when applying the forward F-transform to the initial/boundary condition and the inverse F-transform after the construction of $U_{ij}$. F-transform methodology for non-linear equations, like Burgers, presents the same pros and cons of the FD approach. Care must be taken when using the F-transform methodology to solving equations that form shocks because F-transform methodology is, in its essence, a FD scheme. Besides, non-linear problems violate the continuity assumption. In this case, an alternative (hybrid) approach should be considered to overcome this issue.

The convergence of the F-transform $U_{ij}$ to the wave-equation solution at the grid nodes is assured by the following result.

{\bf Result 1:} Let $u(x,t)$ be a solution of equation (\ref{wave}) and suppose that $u$ is four times continuously differentiable with respect to $x$ and $t$, $r=\displaystyle\frac{c h_t}{h_x}$, and $0<r \leq 1$. Let $u_{ij}=u(x_i,y_j)$. Then
$$\norm{u_{ij}-U_{ij}}=\max_{\substack{i=1,\ldots,n\\ j=1,\ldots,m}}|u_{ij}-U_{ij}|=O(h_t^2)+O(h_x^2),$$
where $u_{ij}=u(x_i,t_j)$ and $U_{ij}$ is the approximative solution given by algorithm (\ref{recur4}) with respect to some fix basis functions.

Besides, using previous result we prove the convergence of inverse fuzzy transform to exact solution of equation (\ref{wave}).

{\bf Result 2:} Assuming the same hypothesis of result 1, then $\displaystyle u^F_{n,m}\to u$  uniformly in $D$ as $n,m \to \infty$, where $u^F_{n,m}$is the inverse F-transform of $u$ with respect to some fix basis functions.

\section{Numerical tests}

In the previous section we applied the F-transform to the wave equation with initial and boundary conditions. We showed that the F-transformed wave equation and the discretized (by a FD scheme) wave equation have, essentially, the same associated difference equation. Which means, the computation of the F-transformed and the FD solutions are exactly the same, except for the calculation of the F-transform of the initial/boundary conditions and the external force, on the other hand, for the FD scheme the initial/boundary conditions and the external force are merely discretized and then used in the difference scheme.

The F-transform of the initial/boundary conditions and the external force are numerical processes, usually calculated by the trapezoidal rule with three points. It uses few points, but by the definition of the basic functions and the F-transform, the integral calculations are performed at $(x_{i-1},x_{i+1})$, $i=1,2,\ldots,n$ (the partition of $[a,b]$), i.e., these calculations can be a bit costly. After the construction of $U_{ij}$, an inverse F-transform is performed to bring the solution  back to the original domain.

The non-homogeneous wave-equation solution based on F-transform, when the external force is noisy, is presented by {{\v S}t{\v e}pni{\v c}ka} and {Val\'a{\v s}ek} \cite{Stepnicka}. In this work we investigate the wave-equation solution behavior, based on the F-transform methodology, for a homogeneous equation with a noise-contaminated initial data.

In order to keep this numerical test simple, consider the homogeneous version of equation (\ref{wave}) with background solution $c=1m/s$:
\begin{equation}\label{condi}
\left\{\begin{array}{l}
  \displaystyle \parcial{^2 u}{x^2} - \parcial{^2 u}{t^2}=0; \\\\
  \displaystyle u(x,0)=\sin(\pi x), \qquad \parcial{u}{t}(x,0)=0; \\\\
   u(0,t)=0, \qquad u(10,t)=0.
\end{array}\right.
\end{equation}

In this example, the domain is $[0,10]\times [0,10]$. The exact solution for this problem is $$\displaystyle u(x,t)=\frac{1}{2}\sin[\pi(x-t)]+\frac{1}{2}\sin[\pi(x+t)].$$ Initially, we evaluate the wave-equation FD solution to be used as reference for comparisons with the noisy numerical solutions. It must be said we are not presenting, as reference, the exact solution, but the FD solution. It is because FD solutions, in general, present some undesirable numerical characteristics like dispersion and dissipation. Of course the chosen grid takes all assumptions into account to minimize/avoid numerical issues, but these effects happen and could blur the conclusions.

Then, we added amplitude and phase random noise to the initial condition (\ref{condi}), more specifically to the $\sin(\pi x)$, and $u_t(x,0)$ is maintained $0$. After that we solved the homogeneous version of (\ref{wave}) with noisy initial conditions using the standard centered FD approximation, the goal with this simulations is to present a solution with no denoise.  A comparison of the F-transform solution with the noise-free FD solution is made to show the robustness of this methodology.

In this numerical test, we uniformly partitioned the space interval $[0,10]$ into $n_x=401$ subintervals $I^x_j=(x_{j-1},x_j)$ of length $\Delta x= 10/(n_x-1)=0.025$. The time interval $[0,10]$ was  uniformly partitioned the space interval $[0,10]$ into $n_t=601$ subintervals $I^t_k=(t_{k-1},t_k)$ of length $\Delta t= 10/(n_t-1)=0.01\overline{6}$, satisfying the stability condition $0<r\leqslant 1$. The basic functions $A_i$ and $B_j$ are symmetrical triangles. Another grid must be defined, the grid we solve the wave equation by finite differences. In fact, this grid is not needed, we just create it to run the reference solutions (noise-free FD solution and noisy FD solution), this grid is characterized as $35\times n_x$ in the space direction and $35\times n_t$ in the time direction.

\begin{figure}[H]
\centerline{ \includegraphics[trim={1.5cm 7.4cm 1.5cm 6.8cm},clip,scale=0.60]{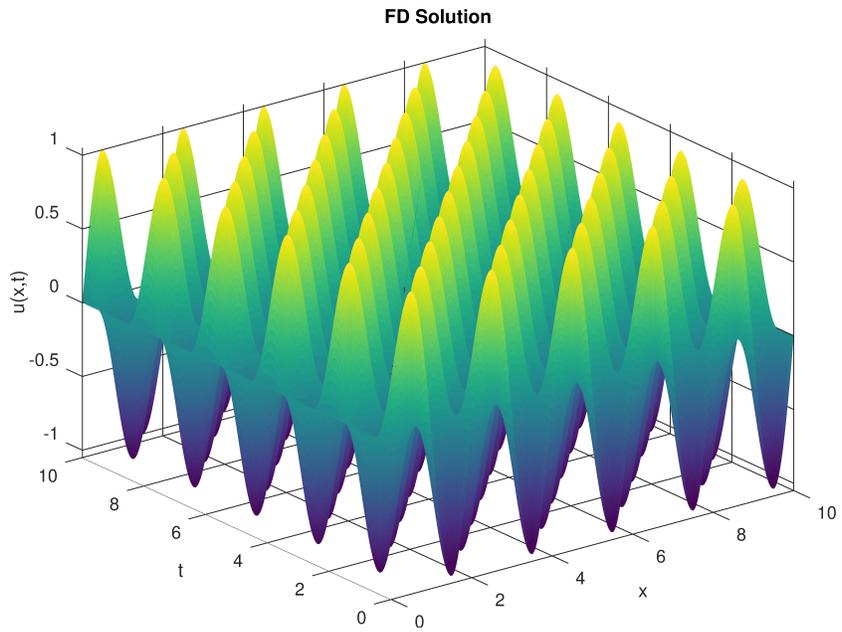}}
\caption{FD solution.}
\label{fig2}
\end{figure}

\begin{figure}[H]
\centerline{ \includegraphics[trim={1.5cm 7.4cm 1.5cm 6.8cm},clip,scale=0.60]{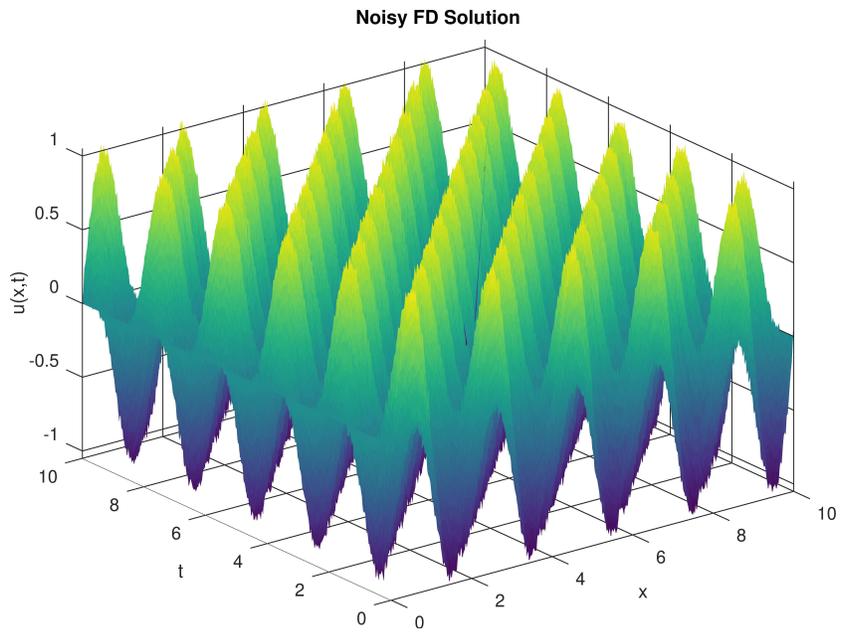}}
\caption{Noisy FD solution.}
\label{fig3}
\end{figure}

\begin{figure}[H]
\centerline{ \includegraphics[trim={1.5cm 7.4cm 1.5cm 6.8cm},clip,scale=0.60]{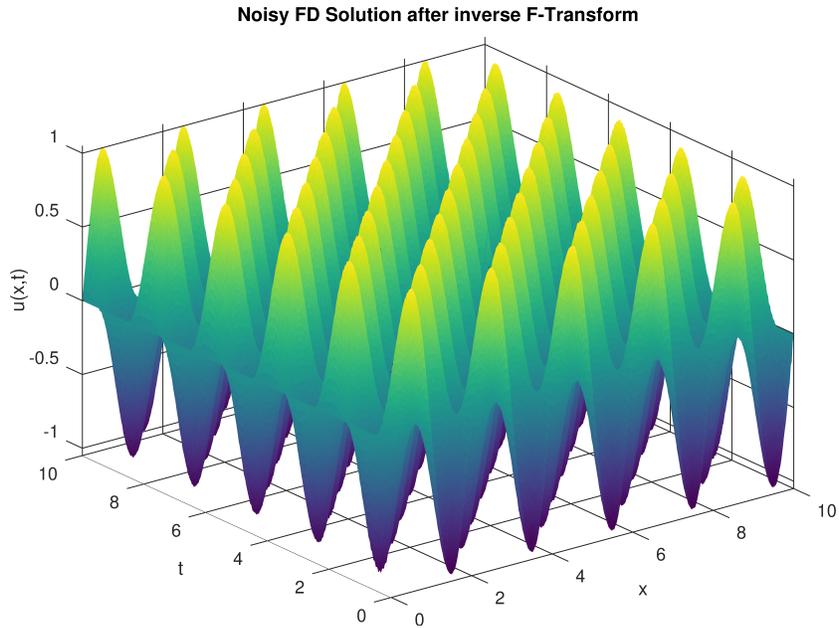}}
\caption{F-transform solution.}
\label{fig4}
\end{figure}

Now, we present the numerical experiment. Figure \ref{fig2} shows the noise-free FD solution surface. In figure \ref{fig3} we present the noisy FD solution, i.e., the FD solution after addition of random noise to the initial condition. The solution surface is a bit crispy. Then, we present, on figure \ref{fig4}, the F-transform solution.

Firstly, we construct $U_{ij}$ using the difference equation (\ref{recur4}) on the $n_x\times n_t$ grid. Secondly, we apply the inverse F-transform to bring the solution back $T-X$ domain.

\begin{figure}[H]
\centerline{ \includegraphics[trim={1cm 1cm 1cm 1cm},clip,scale=0.6]{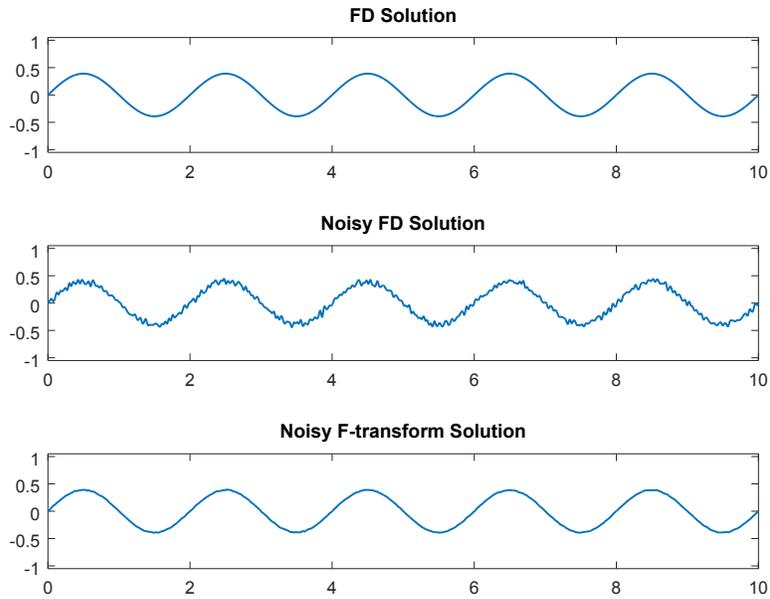}}
\caption{Solution slice at $t=3.63$ s.}
\label{fig7}
\end{figure}

\begin{figure}[H]
\centerline{ \includegraphics[trim={1cm 1cm 1cm 1cm},clip,scale=0.6]{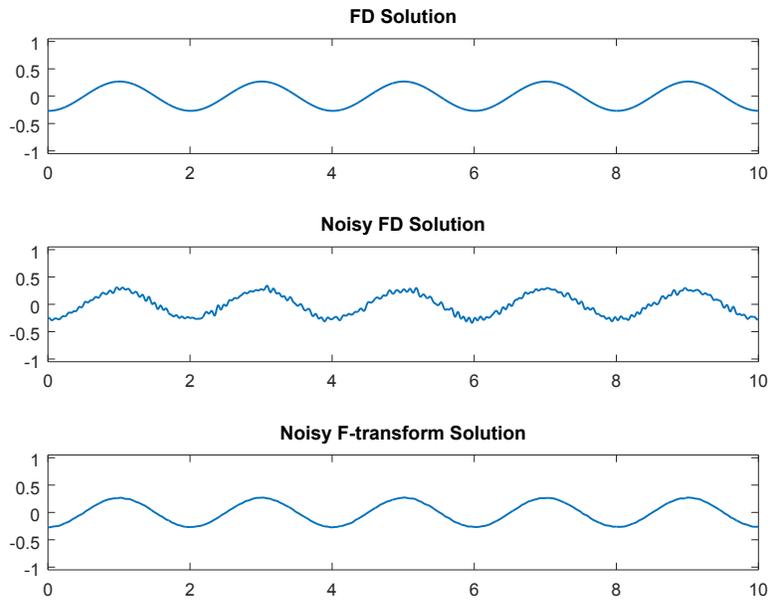}}
\caption{Solution slice at $x=9.08$ m.}
\label{fig8}
\end{figure}

The solution surfaces presented above do not clearly show the effect of the proposed methodology. For this reason, we present on figure \ref{fig7} a time slice at $3.63s$ and on figure \ref{fig8} a space slice at $9.08m$. In these figures we present on top the FD solution, the middle figure is the noisy FD solution, i.e., the solution we would obtain if no other denoise procedure were applied. The bottom graphs present the F-transform solution, which present a good noise attenuation effect.

\begin{figure}[H]
\centerline{ \includegraphics[trim={1cm 1cm 1cm 1cm},clip,scale=0.55]{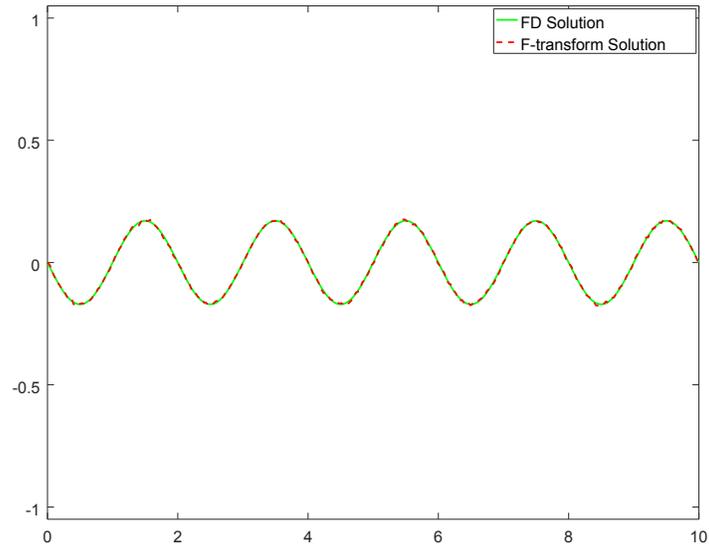}}
\caption{Solution slice at $t=5.45$ s.}
\label{fig9}
\end{figure}

\begin{figure}[H]
\centerline{ \includegraphics[trim={1cm 1cm 1cm 1cm},clip,scale=0.55]{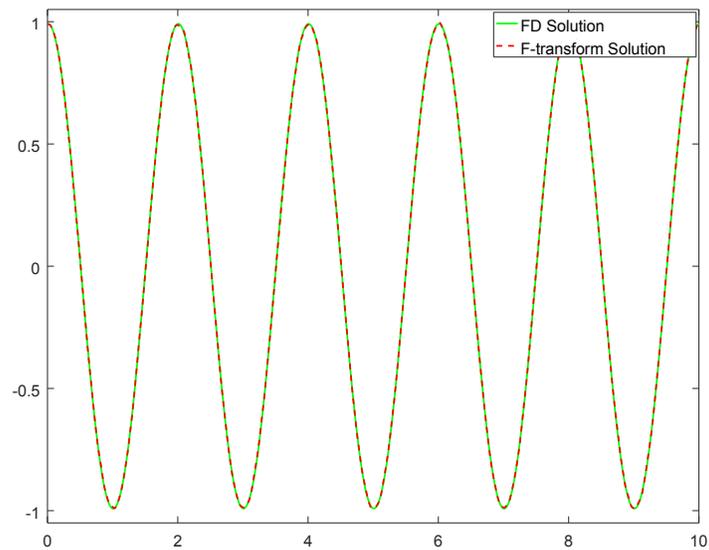}}
\caption{Solution slice at $x=4.54$ m.}
\label{fig10}
\end{figure}

Next pair of slices, figures \ref{fig9} and \ref{fig10}, present an overlay of the noise-free FD solution and the F-transform solution at $5.45s$ and $4.54m$, respectively. The goal of these overlays is to show the effectiveness of the F-transform as a denoise procedure to be adopted during the wave propagation. Amplitude and phase are maintained in the F-transform solution.

The objective of this methodology is to replace the standard FD approximation by a more robust procedure, which is essentially a FD scheme, that attenuates eventual random noise during the wave propagation. Migration and multiple attenuation procedure can get an immediate advantage with the F-transform methodology.

\section{Conclusion}
In this work, we show that F-transform methodology can be used on FD schemes for a certain class of PDEs when the initial condition of such equations is noisy. The F-transform methodology shows noise reduction when compared to the classical FD methods which are, in general, not designed to handle noise.

\bibliographystyle{plain}
\bibliography{paper}
\end{document}